\documentclass[pdflatex,sn-mathphys-num]{sn-jnl}

\usepackage{amsmath,amssymb,amsthm,amsfonts}
\usepackage{mathrsfs}
\usepackage{stmaryrd}
\usepackage{tikz-cd}

\theoremstyle{thmstyleone}
\newtheorem{thm}{Theorem}[section]
\newtheorem{propo}[thm]{Proposition}
\newtheorem{lem}[thm]{Lemma}
\newtheorem{cor}[thm]{Corollary}

\theoremstyle{thmstyletwo}
\newtheorem{ex}[thm]{Example}
\newtheorem{rmk}[thm]{Remark}

\theoremstyle{thmstylethree}
\newtheorem{dfn}[thm]{Definition}

\numberwithin{equation}{section}

\newcommand{\nc}{\newcommand}

\nc{\theo}{\begin{thm}} \nc{\xtheo}{\end{thm}}
\nc{\prop}{\begin{propo}} \nc{\xprop}{\end{propo}}
\nc{\lemm}{\begin{lem}} \nc{\xlemm}{\end{lem}}
\nc{\coro}{\begin{cor}} \nc{\xcoro}{\end{cor}}

\nc{\exam}{\begin{ex}} \nc{\xexam}{\end{ex}}
\nc{\rema}{\begin{rmk}} \nc{\xrema}{\end{rmk}}
\nc{\defi}{\begin{dfn}} \nc{\xdefi}{\end{dfn}}

\nc{\pf}{\begin{proof}} \nc{\xpf}{\end{proof}}

\nc{\on}{\operatorname}
\nc{\fraka}{{\mathfrak a}} \nc{\bba}{{\mathbf a}}
\nc{\frakb}{{\mathfrak b}}
\nc{\frakc}{{\mathfrak c}}
\nc{\frakd}{{\mathfrak d}}
\nc{\frake}{{\mathfrak e}}
\nc{\frakf}{{\mathfrak f}}
\nc{\frakg}{{\mathfrak g}}
\nc{\frakh}{{\mathfrak h}}
\nc{\fraki}{{\mathfrak i}}
\nc{\frakj}{{\mathfrak j}}
\nc{\frakk}{{\mathfrak k}}
\nc{\frakl}{{\mathfrak l}}
\nc{\frakm}{{\mathfrak m}}
\nc{\frakn}{{\mathfrak n}}
\nc{\frako}{{\mathfrak o}}
\nc{\frakp}{{\mathfrak p}}
\nc{\frakq}{{\mathfrak q}}
\nc{\frakr}{{\mathfrak r}}
\nc{\fraks}{{\mathfrak s}}
\nc{\frakt}{{\mathfrak t}}
\nc{\fraku}{{\mathfrak u}}
\nc{\frakv}{{\mathfrak v}}
\nc{\frakw}{{\mathfrak w}}
\nc{\frakx}{{\mathfrak x}}
\nc{\fraky}{{\mathfrak y}}
\nc{\frakz}{{\mathfrak z}}
\nc{\frakA}{{\mathfrak A}}
\nc{\frakB}{{\mathfrak B}}
\nc{\frakC}{{\mathfrak C}}
\nc{\frakD}{{\mathfrak D}}
\nc{\frakE}{{\mathfrak E}}
\nc{\frakF}{{\mathfrak F}}
\nc{\frakG}{{\mathfrak G}}
\nc{\frakH}{{\mathfrak H}}
\nc{\frakI}{{\mathfrak I}}
\nc{\frakJ}{{\mathfrak J}}
\nc{\frakK}{{\mathfrak K}}
\nc{\frakL}{{\mathfrak L}}
\nc{\frakM}{{\mathfrak M}}
\nc{\frakN}{{\mathfrak N}}
\nc{\frakO}{{\mathfrak O}}
\nc{\frakP}{{\mathfrak P}}
\nc{\frakQ}{{\mathfrak Q}}
\nc{\frakR}{{\mathfrak R}}
\nc{\frakS}{{\mathfrak S}}
\nc{\frakT}{{\mathfrak T}}
\nc{\frakU}{{\mathfrak U}}
\nc{\frakV}{{\mathfrak V}}
\nc{\frakW}{{\mathfrak W}}
\nc{\frakX}{{\mathfrak X}}
\nc{\frakY}{{\mathfrak Y}}
\nc{\frakZ}{{\mathfrak Z}}
\nc{\bbA}{{\mathbb A}}
\nc{\bbB}{{\mathbb B}}
\nc{\bbC}{{\mathbb C}}
\nc{\bbD}{{\mathbb D}}
\nc{\bbE}{{\mathbb E}}
\nc{\bbF}{{\mathbb F}} \nc{\bbf}{{\mathbf f}}
\nc{\bbG}{{\mathbb G}}
\nc{\bbH}{{\mathbb H}}
\nc{\bbI}{{\mathbb I}}
\nc{\bbJ}{{\mathbb J}}
\nc{\bbK}{{\mathbb K}}
\nc{\bbL}{{\mathbb L}}
\nc{\bbM}{{\mathbb M}}
\nc{\bbN}{{\mathbb N}}
\nc{\bbO}{{\mathbb O}}
\nc{\bbP}{{\mathbb P}}
\nc{\bbQ}{{\mathbb Q}}
\nc{\bbR}{{\mathbb R}}
\nc{\bbS}{{\mathbb S}}
\nc{\bbT}{{\mathbb T}}
\nc{\bbU}{{\mathbb U}}
\nc{\bbV}{{\mathbb V}}
\nc{\bbW}{{\mathbb W}}
\nc{\bbX}{{\mathbb X}}
\nc{\bbY}{{\mathbb Y}}
\nc{\bbZ}{{\mathbb Z}}
\nc{\calA}{{\mathcal A}}
\nc{\calB}{{\mathcal B}}
\nc{\calC}{{\mathcal C}}
\nc{\calD}{{\mathcal D}}
\nc{\calE}{{\mathcal E}}
\nc{\calF}{{\mathcal F}}
\nc{\calG}{{\mathcal G}}
\nc{\calH}{{\mathcal H}}
\nc{\calI}{{\mathcal I}}
\nc{\calJ}{{\mathcal J}}
\nc{\calK}{{\mathcal K}}
\nc{\calL}{{\mathcal L}}
\nc{\calM}{{\mathcal M}}
\nc{\calN}{{\mathcal N}}
\nc{\calO}{{\mathcal O}}
\nc{\calP}{{\mathcal P}}
\nc{\calQ}{{\mathcal Q}}
\nc{\calR}{{\mathcal R}}
\nc{\calS}{{\mathcal S}}
\nc{\calT}{{\mathcal T}}
\nc{\calU}{{\mathcal U}}
\nc{\calV}{{\mathcal V}}
\nc{\calW}{{\mathcal W}}
\nc{\calX}{{\mathcal X}}
\nc{\calY}{{\mathcal Y}}
\nc{\calZ}{{\mathcal Z}}

\nc{\scrA}{{\mathscr A}}
\nc{\scrE}{{\mathscr E}}
\nc{\scrR}{{\mathscr R}}

\nc{\Bmu}{\mbox{$\raisebox{-0.59ex}{$l$}\hspace{-0.18em}\mu\hspace{-0.88em}\raisebox{-0.98ex}{\scalebox{2}{$\color{white}.$}}\hspace{-0.416em}\raisebox{+0.88ex}{$\color{white}.$}\hspace{0.46em}$}{}}

\nc{\bnu}{{\bar{ \nu}}}

\nc{\olO}{\bar{\calO}}

\nc{\al}{{\alpha}}
\nc{\be}{{\beta}}
\nc{\ga}{{\gamma}} \nc{\Ga}{{\Gamma}}
 \nc{\hGa}{\hat{\Gamma}}
\nc{\ve}{{\varepsilon}}
\nc{\la}{{\lambda}} \nc{\La}{{\Lambda}}
\nc{\om}{\omega} \nc{\Om}{\Omega}
\nc{\sig}{{\sigma}} \nc{\Sig}{{\Sigma}}

\nc{\tnb}{\psi_{\rm tame}}
\nc{\oM}{\overline{{M}}}
\nc{\op}{{\on{op}}}
\nc{\ad}{{\on{ad}}}
\nc{\alg}{{\on{alg}}}
\nc{\Ad}{{\on{Ad}}}
\nc{\Adm}{{\on{Adm}}} \nc{\aff}{{\on{aff}}}
\nc{\Aut}{{\on{Aut}}}
\nc{\Bun}{{\on{Bun}}}
\nc{\cha}{{\on{char}}}
\nc{\der}{{\on{der}}}
\nc{\Der}{{\on{Der}}}
\nc{\diag}{{\on{diag}}}
\nc{\End}{{\on{End}}}
\nc{\Fl}{{\calF\!\ell}}
\nc{\Tr}{{\on{Transp}}}
\nc{\TR}{{\calT\!\calR}}
\nc{\Gal}{{\on{Gal}}}
\nc{\Gr}{{\on{Gr}}}
\nc{\rH}{{\on{H}}}
\nc{\Hom}{{\on{Hom}}}
\nc{\IC}{{\on{IC}}}
\nc{\id}{{\on{id}}}
\nc{\Id}{{\on{Id}}}
\nc{\ind}{{\on{ind}}}
\nc{\Ind}{{\on{Ind}}}
\nc{\Lie}{{\on{Lie}}}
\nc{\Pic}{{\on{Pic}}}
\nc{\pr}{{\on{pr}}}
\nc{\Res}{{\on{Res}}}
\nc{\res}{{\on{res}}} \nc{\Sat}{{\on{Sat}}}
\nc{\s}{{\on{sc}}}
\nc{\drv}{{\on{der}}}
\nc{\sgn}{{\on{sgn}}}
\nc{\Spec}{{\on{Spec}}}\nc{\Spf}{\on{Spf}}
\nc{\Sph}{\on{Sph}}
\nc{\St}{{\on{St}}}
\nc{\tr}{{\on{tr}}}
\nc{\Mod}{{\mathrm{-Mod}}}
\nc{\Hilb}{{\on{Hilb}}}
\nc{\Ext}{{\on{Ext}}}
\nc{\vs}{{\on{Vec}}}
\nc{\ev}{{\on{ev}}}
\nc{\nO}{{\breve{\calO}}}
\nc{\tS}{{\tilde{S}}}
\nc{\spe}{{\on{sp}}}
\nc{\loc}{{\on{loc}}}
\nc{\Sym}{{\on{Sym}}}
\nc{\Cone}{{\on{C}}}
\nc{\syn}{{\on{syn}}}
\nc{\reg}{{\on{reg}}}
\nc{\colim}{{\on{colim}}}
\nc{\Norm}{{\on{N}}}

\nc{\nscrR}{{\mathscr{R}^{\on{nr}}}}

\nc{\GL}{{\on{GL}}}
\nc{\U}{{\on{U}}}
\nc{\Gl}{\on{Gl}}
\nc{\GSp}{{\on{GSp}}}
\nc{\gl}{{\frakg\frakl}}
\nc{\SL}{{\on{SL}}}
\nc{\SU}{{\on{SU}}}
\nc{\SO}{{\on{SO}}}
\nc{\PGL}{{\on{PGL}}}

\nc{\Conv}{{\on{Conv}}}
\nc{\Rep}{{\on{Rep}}}
\nc{\Dom}{{\on{Dom}}}
\nc{\red}{{\on{red}}}
\nc{\act}{{\on{act}}}
\nc{\nr}{{\on{nr}}}
\nc{\ctf}{{\on{ctf}}}

\nc{\str}{{\on{-}}}
\nc{\os}{{\bar{s}}}
\nc{\oeta}{{\bar{\eta}}}

\nc{\hookto}{\hookrightarrow}
\nc{\longto}{\longrightarrow}
\nc{\leftto}{\leftarrow}
\nc{\onto}{\twoheadrightarrow}
\nc{\lonto}{\twoheadleftarrow}

\nc{\uG}{{\underline{G}}}
\nc{\uA}{{\underline{A}}}
\nc{\uS}{{\underline{S}}}
\nc{\uT}{{\underline{T}}}
\nc{\uM}{{\underline{M}}}
\nc{\uP}{{\underline{P}}}
\nc{\uB}{{\underline{B}}}
\nc{\uN}{{\underline{N}}}

\nc{\ucG}{{\underline{\calG}}}
\nc{\ucA}{{\underline{\calA}}}
\nc{\ucS}{{\underline{\calS}}}
\nc{\ucT}{{\underline{\calT}}}
\nc{\ucalM}{{\underline{\calM}}}
\nc{\ucP}{{\underline{\calP}}}
\nc{\ucalN}{{\underline{\calN}}}

\nc{\bF}{{\breve{F}}}

\nc{\oFl}{{\overline{\Fl}}}
\nc{\bU}{{\overline{U}}}
\nc{\tGr}{{\tilde{\Gr}}}
\nc{\cGr}{\calG\! r}
\nc{\oGr}{\overline{\on{Gr}}}
\nc{\ocGr}{\overline{\calG\! r}}
\nc{\co}{{\colon}}
\nc{\sch}[1]{(Sch/{#1})}
\nc{\HypLoc}[1]{HypLoc({#1})}

\nc{\ohtimes}{\stackrel{!}{\otimes}}
\nc{\boxtilde}{\widetilde{\boxtimes}}
\nc{\vstar}{{\varhexstar}}

\nc{\Div}{\on{Div}}
\nc{\Sht}{\on{Sht}}
\nc{\Frob}{\on{Frob}}

\nc{\x}{\times}
\nc{\bsl}{\backslash}
\nc{\algQl}{{\bar{\bbQ}_\ell}}
\nc{\sF}{{\bar{F}}}
\nc{\nF}{{\breve{F}}}
\nc{\nW}{{W^{\on{nr}}}}
\nc{\sk}{{\bar{k}}}
\nc{\cont}{\on{c}}
\nc{\Supp}{\on{Supp}}
\nc{\blt}{\bullet}
\nc{\dom}{\on{dom}}
\nc{\scon}{{\on{sc}}}
\nc{\Affine}{\on{Aff}}
\nc{\nscrA}{\mathscr{A}^{\on{nr}}}
\nc{\nfraka}{{\bbf^{\on{nr}}}}
\nc{\ran}{{\rangle}}
\nc{\lan}{{\langle}}
\nc{\bk}{{\bar{k}}}
\nc{\tF}{{\tilde{F}}}
\nc{\sS}{{\bar{S}}}
\nc{\LG}{{^\text{L}\hspace{-0.04cm}G}}
\nc{\LL}{{^\text{L}\hspace{-0.07cm}L}}
\nc{\et}{{\text{\rm \'et}}}
\nc{\inv}{{\on{inv}}}
\nc{\Hecke}{{\on{Hecke}}}
\nc{\Isom}{{\on{Isom}}}
\nc{\oSht}{{\overline{\on{Sht}}}}
\nc{\umu}{{\underline \mu}}
\nc{\AIJ}{{\calO_X[{\scriptstyle{\calI\over \calJ}}]}}
\nc{\Proj}{{\on{Proj}}}
\nc{\Bl}{{\on{Bl}}}

\nc{\Pos}{{\on{Pos}}}
\nc{\Sets}{{\on{Sets}}}
\nc{\AffSch}{{\on{AffSch}}}
\nc{\Groups}{{\on{Groups}}}
\nc{\Gpds}{{\on{Groupoids}}}
\nc{\Sch}{{\on{Sch}}}
\nc{\fl}{{\on{flat}}}

\nc{\pot}[1]{ [\hspace{-0,5mm}[ {#1} ]\hspace{-0,5mm}] }
\nc{\rpot}[1]{ (\hspace{-0,7mm}( {#1} )\hspace{-0,7mm}) }

\nc{\defined}{\hspace{0.1cm}\stackrel{\text{\tiny \rm def}}{=}\hspace{0.1cm}}

\begin{document}

\title[Algebraic Magnetism: $X$-products of attractors]{Algebraic Magnetism: $X$-products of attractors via $\mathbb{F}_1$-geometry}

\author[1]{\fnm{Arnaud} \sur{Mayeux}}\email{mayeux@wisc.edu}

\affil[1]{\orgname{University of Wisconsin--Madison}, \orgaddress{\city{Madison}, \state{WI}, \country{USA}}}

\abstract{We give a formula for $X$-products of attractors using $\mathbb{F}_1$-geometry.}



\maketitle

\section{Introduction}

Tits heuristic idea of the field with one element $\mathbb{F}_1$ and Deitmar’s monoid-based $\mathbb{F}_1$-geometry both aim to describe aspects of algebraic geometry through underlying combinatorial structures.
Algebraic magnetism follows a similar philosophy by encoding geometric dynamics via commutative monoids.
In this note, we answer a purely magnetic question in terms of $\mathbb{F}_1$-geometry. Algebraic magnetism was announced in \cite{Ma0} and developed in \cite{Ma}; see also \cite{Ma1}. This note builds on and complements \cite{Ma}. Our starting point is the following question: given submonoids of an abelian group $N_1, \ldots, N_d \subset M$, viewed as magnets for an algebraic action of $D(M)$ on an algebraic space $X$, can we understand and describe
$
X^{N_1} \times_X \cdots \times_X X^{N_d}\,?
$
We provide an answer to this question in the case where $N_1, \ldots, N_d$ satisfy $N_1^{\mathrm{gp}} = \cdots = N_d^{\mathrm{gp}}$. In this situation, the $N_i$ together give rise to a scheme over $\mathbb{F}_1$, denoted $\underline{\mathcal{N}}$ (Proposition \ref{prop:schoverfun}), and to a scheme $A(\mathcal{N})$ (Propositions \ref{prop:pushsch} and \ref{prop:pushequiv}) endowed with an action of $D(M)$, such that
\[
X^{N_1} \times_X \cdots \times_X X^{N_d} \cong \underline{\mathrm{Hom}}^{D(M)}\bigl(A(\mathcal{N}), X\bigr) ~~~~~~\text{ (Theorem \ref{thm:fiberprodmono})}.
\]
We then show that the space $X^{\mathcal{N}}:= \underline{\mathrm{Hom}}^{D(M)}\bigl(A(\mathcal{N}), X\bigr)$ is representable under natural assumptions and that, for $X$ separated, the morphism $X^{\mathcal{N}}\to X$ is a monomorphism. We also provide examples and relate it to existing constructions in the literature (Example \ref{examm}).

\section{Schemes over $\bbF _1$}
Schemes over $\mathbb{F}_1$ \cite{Dei}, also called monoschemes \cite{Og18}, are sheaves of monoids on topological spaces that are built locally from commutative monoids. In the following, we recall some main structural definitions of these objects, and refer to \cite{Dei} and \cite{Og18} for details.
\defi
If  $Q$ is a monoid, then
$a_Q := \mathrm{spec}(Q)$
is the (locally) monoidal space consisting of the set of prime ideals of
$Q$, endowed with the Zariski topology and the sheaf of monoids defined in \cite{Og18}.
\xdefi

\defi A monoidal space is a pair $(X, \calM _X)$, where $X$ is a topological
space and $\calM _X$ is a sheaf of commutative monoids on $X.$ A morphism of
monoidal spaces
$(f, f^\sharp ): (X,\calM_X) \to (Y, \calM_Y )
$
is a pair $(f, f^\sharp )$, where $f : X \to Y$ is a continuous function and
$f ^\sharp : f^{-1} (\calM _Y ) \to \calM _X$
is a homomorphism of sheaves of monoids such that for each point $x \in X$,
the stalk $f_x^\sharp : \calM _{Y , f (x)} \to  \calM_{X, x}$ of $f^\sharp$ at $x$ is a local homomorphism of monoids.
\xdefi

\defi A locally monoidal space is affine if it is isomorphic to
$\mathrm{spec}(Q)$ for some monoid $Q.$
A monoscheme (or a scheme over $\mathbb{F}_1$) is a locally monoidal space that
admits an open cover by affines, and the category of monoschemes is the full
subcategory of the category of locally monoidal spaces consisting of such
spaces.
\xdefi
Given a scheme over $\mathbb{F}_1$, say $X$, one can associate to it a scheme over $\mathbb{Z}$ (cf. \cite{Dei}). This construction is denoted by
$
X \otimes_{\mathbb{F}_1} \mathbb{Z},
$
and may be viewed as an extension of scalars from $\mathbb{F}_1$ to $\mathbb{Z}$.
On affines, the construction is particularly simple. If
$
X=\operatorname{spec}( Q)
$
for a commutative monoid $Q$, then
$
X\otimes_{\mathbb{F}_1}\mathbb{Z}
 \;=\;
\operatorname{Spec}(\mathbb{Z}[Q]),
$
where $\mathbb{Z}[Q]$ denotes the monoid algebra of $Q$ over $\mathbb{Z}$. The general construction is obtained by gluing these affine pieces.

\section{Monoids, $\bbF _1$-schemes and schemes}
If $N$ is a commutative monoid, recall that $A(N)=\Spec( \bbZ [N])$
denotes the diagonalizable monoid scheme over $\bbZ$ attached to it \cite{Ma}. If $N$ is moreover a group, $A(N)$ is endowed with a group scheme structure denoted $D(N)$.

Let $M$ be a finitely generated abelian group.
Let $\calN =\{N_{\tau}\}_{\tau \in \mathscr{A}}$ be a set of fine submonoids of $M$ such that $N_{\tau}^{gp} = N_{\tau '} ^{gp}=: \calN ^{gp}$ for all $\tau, \tau' \in \mathscr{A}$. Note that for all $\tau \in \mathscr{A}$, $D(\calN ^{gp})\to A (N_{\tau})$ is an open immersion (cf. e.g. \cite[p.~67]{Og18}).

\prop \label{prop:pushsch}
The diagram
\[
D(\mathcal N^{{gp}})\longrightarrow A(N_\tau), \qquad \tau\in\mathscr A,
\]
admits a colimit in the category of schemes; this colimit is denoted by $A(\mathcal N)$.
\xprop

\pf
Since each $D(\mathcal N^{\mathrm{gp}})\to A(N_\tau)$ is an open immersion, the schemes $A(N_\tau)$ glue along the common open subscheme $D(\mathcal N^{\mathrm{gp}})$; by the standard universal property of gluing along open immersions, this glued scheme is the colimit of the diagram.
\xpf

\prop \label{prop:schoverfun}There exists a scheme over $\mathbb{F} _1$, denoted $\underline{\calN}$ such that $\underline{\calN} \otimes _{\mathbb{F} _1} \bbZ \cong A (\calN)$.
\xprop
\pf
Define $\underline{\mathcal N}$ by gluing the affine monoschemes $\operatorname{spec}(N_\tau)$ along the common open subspace $\operatorname{spec}(\mathcal N^{\mathrm{gp}})$. This is well-defined since all $N_\tau$ have the same group completion. Hence $\underline{\mathcal N}$ is a monoscheme, i.e. a scheme over $\mathbb F_1$.
Applying base change $N \mapsto \mathbb Z[N]$ (cf. \cite{Dei}), localization of monoids commutes with monoid algebras, so the gluing diagram becomes the one defining $A(\mathcal N)$. Therefore
$
\underline{\mathcal N} \otimes_{\mathbb F_1} \mathbb Z \cong A(\mathcal N).
$
\xpf
We now put $Z= \calN ^{gp} \subset M$. For each $\tau$, the morphism $D(Z) \to A(N_\tau)$ is $D(Z)$-equivariant.  We thus obtain an action of $D(Z)$ on $A(\calN)$.

\prop \label{prop:pushequiv}
The diagram
\[
D(\mathcal N^{{gp}})\longrightarrow A(N_\tau), \qquad \tau\in\mathscr A,
\]
admits a colimit in the category of schemes with a $D(M)$-action. This colimit is isomorphic to $A(\mathcal N)$ as scheme.
\xprop

\pf
Each $A(N_\tau)$ carries a natural action of $D(M)$ (since $N_\tau \subset M$), and the morphisms
$
D(Z)\to A(N_\tau)
$
are $D(M)$-equivariant.
The gluings along $D(Z)$ are compatible with the $D(M)$-actions, hence $A(\mathcal N)$ inherits a $D(M)$-action such that all maps
$
A(N_\tau)\to A(\mathcal N)
$
are $D(M)$-equivariant.
Now let $X$ be a $D(M)$-scheme and suppose given $D(M)$-equivariant morphisms
$
f_\tau : A(N_\tau)\to X
$
compatible on $D(Z)$. By the universal property in schemes, they glue to a unique morphism
$
f : A(\mathcal N)\to X.
$
Since the $f_\tau$ are $D(M)$-equivariant and the action on $A(\mathcal N)$ is defined by gluing, it follows that $f$ is $D(M)$-equivariant.
Thus $A(\mathcal N)$ satisfies the required universal property in the category of $D(M)$-schemes.
\xpf

\section{$X$-product of attractors}
We proceed with the setting of the previous section.
We work over $S= \Spec ( \bbZ)$. We use the notation of \cite{Ma}, in particular, in the following, the notation $X^N $ means $\underline{\Hom} ^{D(M)} (A(N),X)$ (i.e. the functor $(T \to S) \mapsto \Hom _T ^{D(M)_T} (A(N)_T,X_T)$, where $\Hom _T ^{D(M)_T}$ denotes the set of $D(M)_T$-equivariant morphisms of $T$-spaces).
\theo \label{thm:fiberprodmono} Let $X$ be an algebraic space endowed with an action of $D(M)$. Assume that one of the following conditions holds
\begin{enumerate}
\item $X$ is quasi-separated and locally of finite presentation over $S$ and $D(M/Z)$ has connected fibers (e.g. $Z=M$), or
\item $X$ is separated.
\end{enumerate}
Then $${\prod _{\tau \in \mathscr{A} }}{}_X X^{N_\tau}  =... \times _X X^{N_\tau} \times _X X^{N_{\tau '}}\times _X ...= \underline{\Hom}^{D(M)} (A(\calN) , X).$$
\xtheo
\pf
Let $T \to S$ be a test scheme. For each $\tau$, let $g_\tau : D(M)\to A(N_\tau)$ be the morphism of diagonalizable monoid schemes dual to $N_\tau \subset M$, and let $h_\tau : D(Z) \to  A(N_\tau )$ be the morphism of diagonalizable monoid schemes dual to $N_\tau \subset Z$. Similarly, let $k : D(M) \to D(Z)$ be the morphism dual to $Z \subset M$. Since $X^M\cong X$ (cf. \cite{Ma}), we have:
\begin{align*}
&\Hom  (T , ... \times _X X^{N_\tau} \times _X X^{N_{\tau '}}\times _X ...)
\cong   \Hom  (T , ... \times _{X^M} X^{N_\tau} \times _{X^M} X^{N_{\tau '}}\times _{X^M} ...) \\
\cong &  \{ \{f_{\tau } \in \Hom ^{D(M)_T} ( A(N_\tau)_T , X_T)  \}_{\tau \in \mathscr{A}} | f_{\tau} \circ g_\tau = f_{\tau'} \circ g_{\tau '} ~\forall \tau , \tau' \in \mathscr{A} \} \\
\cong &  \{ \{f_{\tau } \in \Hom ^{D(M)_T} ( A(N_\tau)_T , X_T)  \}_{\tau \in \mathscr{A}} | f_{\tau} \circ h_\tau \circ k = f_{\tau'} \circ h_{\tau '} \circ k ~\forall \tau , \tau' \in \mathscr{A} \}.
\end{align*}
We now prove that we have
$$
   \{ \{f_{\tau } \in \Hom ^{D(M)_T} ( A(N_\tau)_T , X_T)  \}_{\tau \in \mathscr{A}} | f_{\tau} \circ h_\tau \circ k = f_{\tau'} \circ h_{\tau '} \circ k ~\forall \tau , \tau' \in \mathscr{A} \}
$$
$$
\cong  \{ \{f_{\tau } \in \Hom ^{D(M)_T} ( A(N_\tau)_T , X_T)  \}_{\tau \in  \mathscr{A}} | f_{\tau} \circ h_\tau = f_{\tau'} \circ h_{\tau '}  ~\forall \tau , \tau' \in \mathscr{A} \} .$$
It is enough to prove that $X^Z \to X=X^M$ is a monomorphism.
If $X$ is separated over $S$, this is proved in \cite{Ma}.
If $X$ is quasi-separated over $S$ and $D(M/Z)$ has connected fibers, it is enough to observe that $X^Z =X^{D(M/Z)}$ (cf. \cite[Proposition 3.35]{Ma}) and apply \cite[Theorem A.1]{Ma}.
Now Proposition \ref{prop:pushequiv} implies that $$\{ \{f_{\tau } \in \Hom ^{D(M)_T} ( A(N_\tau)_T , X_T)  \}_{\tau \in \mathscr{A}} | f_{\tau} \circ h_\tau = f_{\tau'} \circ h_{\tau '}  ~\forall \tau , \tau' \in \mathscr{A} \} \cong {\Hom}^{D(M)_T} (A(\calN) _T, X_T).$$ This finishes the proof.
\xpf

\defi The functor $\underline{\Hom}^{D(M)} (A(\calN) , X)$ is denoted $X^\calN$ and is called the attractor associated to $\underline{\calN}$.
\xdefi

\prop
Assume that one of the following conditions holds
\begin{enumerate}
\item $X \to S$ is affine, or
\item the set $\mathscr{A}$ is finite, each $X^{N_\tau}$ is representable, and one of the following conditions holds
\begin{enumerate}
\item $X$ is quasi-separated and locally of finite presentation over $S$ and $D(M/Z)$ has connected fibers (e.g. $Z=M$), or
\item $X$ is separated.
\end{enumerate}
\end{enumerate}
then $X^\calN$ is representable by an algebraic space over $S$.

Moreover, if $X\to S$ is affine, then $X^\calN = X^{\cap _{\tau \in  \mathscr{A}} N_\tau}$.
\xprop

\pf
If $X \to S$ is affine, then each $X^{N_\tau}$ is representable by a closed subscheme $V(J_\tau)$ of $X=\Spec (A)$ where the ideal $J_\tau \subset A$ is given explicitly in \cite{Ma}, namely $J_\tau = \langle A_m |m \in M \setminus N _\tau\rangle$. This explicit description of $J_\tau$ implies that $X^\calN$ $(\cong \Spec(\otimes _{\tau \in \mathscr{A} ~A}~ A/J_\tau )\cong \Spec (A/ (\sum _{\tau \in \mathscr{A} } J_\tau ) ))$ equals $ X^{\cap _{\tau \in  \mathscr{A}} N_\tau}$ and is representable.
In the second case, by \cite{Ma}, each $X^{N_\tau}$ is representable. Since $\mathscr{A}$ is finite, Theorem \ref{thm:fiberprodmono} implies that $X^\calN$ is representable.
\xpf

\exam \label{examm}
The formula $X^\calN = X^{\cap _{\tau \in  \mathscr{A}} N_\tau}$ does not extend to the non-affine case.
Let $M=Z=\bbZ$ and $\calN =\{N_{\tau}= \tau \bbN \mid \tau \in \{+,-\} \}$.
 Then $A(\calN ) \cong  \bbP ^1$ with the action $t\cdot [x:y]=[tx:y]$. Let us take $X= A(\calN )$.
 By \cite[Example 1.6.4]{DG}, $X^{+\bbN}\to X$ is $(\bbP^1\setminus\{\infty\})\sqcup\{\infty\}\to\bbP^1$ and, symmetrically, $X^{-\bbN}\to X$ is $(\bbP^1\setminus\{0\})\sqcup\{0\}\to\bbP^1$. Since fiber products commute with disjoint unions and fiber products of immersions are intersections, Theorem \ref{thm:fiberprodmono} gives
 $$X^{\calN} \cong \bbG_m\sqcup\{0\}\sqcup\{\infty\} ~~\not \cong  ~~ X^{-\bbN  \cap + \bbN } = X^{0} \cong \{0, \infty\},$$
 the $\bbG_m$-equivariant maps $\bbP^1\to\bbP^1$ being $z\mapsto az$ ($a$ a unit) and the two constant maps.
Examples of attractors of the form $Y^\calN$ where $\calN = \{N_{\tau}= \tau \bbN |\tau \in \{+,-\} \}$ appear in \cite[§1.9, Proposition 1.9.4, Remark 1.9.5]{DG}.
\xexam

\prop Assume that $\mathscr{A}$ is finite.
If $X/S$ is separated, then $X^\calN \to X$ is a monomorphism.
\xprop
\pf
We prove it by induction on $\#\mathscr{A}$. If $\#\mathscr{A} =1$, this is proved in \cite{Ma}. Otherwise we have $\mathscr{A} = \mathscr{A'}\sqcup \{\tau\} $ with $\tau \in \mathscr{A}$.
We have $X^{\calN } \cong  X^{\calN'} \times _X X^{N_\tau}$. By induction, $X^{\calN'} \to X$ is a monomorphism, so by \cite[\href{https://stacks.math.columbia.edu/tag/042P}{Tag 042P}]{stacks-project}, $X^{\calN } \cong  X^{\calN'} \times _X X^{N_\tau} \to X^{N_\tau}$ is a monomorphism. By \cite{Ma}, $X^{N_\tau} \to X$ is a monomorphism. So by  \cite[\href{https://stacks.math.columbia.edu/tag/042O}{Tag 042O}]{stacks-project}, $X^\calN \to X$ is a monomorphism.
\xpf

\bmhead{Acknowledgements}
The author thanks the referee for a careful reading of the manuscript and for valuable suggestions that improved both the content and the presentation of the paper.


\begin{thebibliography}{9}
\small
\bibitem{Dei} Deitmar, A.: Schemes over $\mathbb{F}_1$. In:  van der Geer, G., Moonen, B., Schoof, R. (eds.) Number Fields and Function Fields---Two Parallel Worlds. Progress in Mathematics, vol. 239, pp. 87--100. Birkh\"auser, Boston (2005) \href{https://doi.org/10.1007/0-8176-4447-4_6}{doi:10.1007/0-8176-4447-4\_6}

\bibitem{DG} Drinfeld, V., Gaitsgory, D.: On a theorem of Braden. Transform. Groups \textbf{19}(2), 313--358 (2014) \href{https://doi.org/10.1007/s00031-014-9267-8}{doi:10.1007/s00031-014-9267-8}

\bibitem{Ma0} Mayeux, A.: Magnets and attractors of diagonalizable group schemes actions. In: Dobrev, V. (ed.) Lie Theory and Its Applications in Physics. Springer Proc. Math. Stat., vol. 473, pp. 437--442. Springer, Singapore (2025) \href{https://doi.org/10.1007/978-981-97-6453-2_36}{doi:10.1007/978-981-97-6453-2\_36}

\bibitem{Ma} Mayeux, A.: Algebraic magnetism. Semigroup Forum \textbf{113}(1), 126--201 (2026) \href{https://doi.org/10.1007/s00233-026-10645-3}{doi:10.1007/s00233-026-10645-3}

\bibitem{Ma1} Mayeux, A.: Algebraic magnetism invariants of self-actions of diagonalizable monoid schemes. Semigroup Forum \textbf{112}(2), 532--540 (2026) \href{https://doi.org/10.1007/s00233-025-10572-9}{doi:10.1007/s00233-025-10572-9}

\bibitem{Og18} Ogus, A.: Lectures on Logarithmic Algebraic Geometry. Cambridge Studies in Advanced Mathematics, vol. 178. Cambridge University Press, Cambridge (2018) \href{https://doi.org/10.1017/9781316941614}{doi:10.1017/9781316941614}

\bibitem{stacks-project} The Stacks Project authors: The Stacks Project. \url{https://stacks.math.columbia.edu}

\end{thebibliography}
\end{document}